\input amstex
\input epsf
\documentstyle{amsppt}
\NoBlackBoxes
\nologo
\tolerance 1000
\vsize=227mm
\hsize=152mm

\def\P{{\Cal P}}
\def\R{{\Cal R}}

\def \e {\varepsilon}
\define \Ga {\Gamma}
\define \la {\lambda}
\define \pa {\partial}
\define \bR {\Bbb R}
\define \bC {\Bbb C}
\define \bZ {\Bbb F}
\define \bF {\Bbb F}
\define \al {\alpha}
\define \siij {\sigma_{i,j}}
\define \Si {\Sigma}
\define \si {\sigma}

\def\fy{\varphi}
\def\ls{\leqslant}
\def\gs{\geqslant}

\topmatter

\title
Periodic de Bruijn triangles: exact and asymptotic results
\endtitle

\author
B.~Shapiro$^{\ddag}$, M.~Shapiro$^*$ and A.~Vainshtein$^{\dag}$
\endauthor

\affil
$^\ddag$ Department of Mathematics, University of Stockholm\\
S-10691, Sweden, {\tt shapiro\@matematik.su.se}\\
$^*$ Department of Mathematics, Michigan State University\\
East Lansing, MI 48824, USA, {\tt mshapiro\@math.msu.edu}\\
$^\dag$ Dept. of Mathematics and Computer Science, University of Haifa
\\ Mount Carmel, 31905 Haifa, Israel, {\tt alek\@mathcs11.haifa.ac.il}
\endaffil

\abstract
We study the distribution of the number of permutations with a given
periodic up-down
sequence w.r.t. the last entry, find exponential generating functions
and prove asymptotic formulas
for this distribution.
\endabstract
\rightheadtext{Periodic de Bruijn triangles}
\leftheadtext{B.~Shapiro, M.~Shapiro and A.~Vainshtein}
\endtopmatter

\document
\heading {\S 1. Introduction and results} \endheading

Let $\si=(\si_1,\dots,\si_n)$ be a permutation of length $n$.
We associate with $\si$ its
{\it up-down sequence\/} (sometimes called the {\it shape\/} of $\si$, or the
{\it signature\/} of $\si$)  $\P(\si)=(p_{1},\ldots,p_{n-1})$, which is a
binary vector of length $n-1$ such that $p_{i}=1$ if $\si_{i}<\si_{i+1}$
and $p_{i}=0$ otherwise. During the last 120 years, many authors have studied
the number $\sharp_n^{\P}$  of all permutations of length $n$ with a given
up-down sequence $\P$. Apparently, for the first time this problem was
investigated by D.~Andr\`e \cite{An1, An2}, who considered the so-called
alternating (or up-down) permutations  corresponding to the
sequence $\P=(1,0,1,0,\ldots )=(10)^*$ and proved that the exponential generating
function for the number of such permutations is equal to $\tan x+\sec x$.
In \cite{An3} he proved that this number grows asymptotically as
$2n!(2/\pi)^{n+1}$.

A general approach to this problem was
suggested by MacMahon (see \cite{MM}). This approach leads to determinantal
formulas for $\sharp_n^{\P}$, rediscovered later by Niven \cite{Ni} from very basic
combinatorial considerations. For the relations of this approach to the
representation theory of the symmetric group, and for
its generalizations, see \cite{Fo, St1, BW}.

Another, purely combinatorial approach to the same problem was suggested by
Carlitz \cite{Ca1}. His general recursive formula for $\sharp_n^{\P}$ is rather
difficult to use. However,  he managed to obtain
explicit expressions for the corresponding exponential generating functions
for certain {\it periodic\/} cases, that is for up-down sequences of the form
$\P=(p)^*$, where $p$ is a binary vector of a fixed length called the {\it period\/}
of $\P$.
In \cite{Ca1, Ca2} he considered the case $\P=(1^k0)^*$ and expressed the corresponding
generating function via the {\it Olivier functions of the $k$th order}
$$
\fy_{k,i}(x)=\sum_{j=0}^\infty \frac{x^{jk+i}}{(jk+i)!}, \qquad 0\ls i\ls k-1.
$$
Another case, $\P=(1^20^2)^*$, was considered in \cite{CS1,CS2} and solved via
Olivier functions of the fourth order.
It follows that asymptotically $\sharp_n^{\P}$ in this case grows as
$4n!(2/\gamma)^{n+1}$,
where $\gamma=3.7502\ldots$ is the smallest positive solution of the equation
$\cos t\cosh t+1=0$.

The general periodic problem was solved completely in \cite{CGJN}. As in the
two particular
cases mentioned above, the answer is expressed via Olivier functions. The techniques
used involves matrix Riccati equations, and is rather complicated. For a
different solution based on M\"obius functions see \cite{St2, Ch.~3.16, and
Ex.~3.80}.

An additional dimension in the problem was introduced by Entringer \cite{En} who
studied the distribution of the alternating permutations by the last entry. He
observed that the number $\sharp_{i,j}$ of alternating permutations of length $i$
whose last entry equals $j$ satisfy the following equations:
$$\alignedat 3
&\sharp_{i,j}=\sharp_{i,j-1}+\sharp_{i-1,j-1}, &\quad &\sharp_{i,1}=0,&
\qquad &i=2k, k>0,\\
&\sharp_{i,j}=\sharp_{i,j+1}+\sharp_{i-1,j}, &\quad &\sharp_{i,i}=0,&
\qquad &i=2k+1, k>0,
\endalignedat\tag1.1
$$
with $\sharp_{1,1}=1$. These equations can be represented
graphically as the following triangle

\vskip 15pt
\centerline{\hbox{\epsfysize=2.3cm\epsfbox{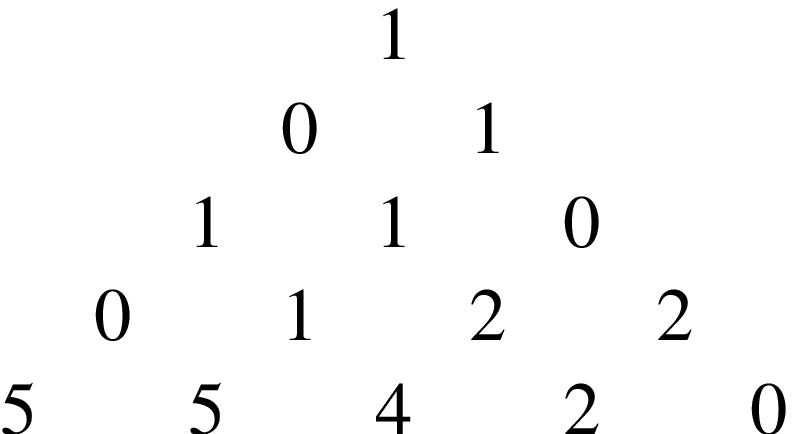}}}
\midspace{0.1mm} \caption{Fig.~1. The Entringer triangle}

Each even row of the triangle starts with $0$, and an entry in such a row is equal to
the sum of its {\it left\/} neighbors in the current and in the previous rows. Similarly,
each odd row (except for the first one) ends with $0$, and an entry in such a
row is equal to
the sum of its {\it right\/} neighbors in the current and in the previous rows.

The Entringer triangle was studied by many authors. In particular, Arnold
\cite{Ar1, Ar2} gave an
interpretation of the entries of this triangle in terms of real polynomials with
real critical values.
Besides, he considered the exponential generating function
$$
A(x,y)=\sum_{i\gs1}\sum_{j=1}^i (-1)^{(i-1)(i-2)/2}\sharp_{i,j}
\frac{x^{i-j}y^{j-1}}{(i-j)!(j-1)!}
$$
and proved that $A(x,y)=e^y/\cosh(x+y)$. In fact, $A(x,y)$ is the generating
function of the
{\it signed Entringer triangle\/}, which is obtained from the ordinary one by
reversing signs
in each $i$th row, where $i$ equals $0$ or $3$ modulo $4$. Observe that the entries
$\tilde\sharp_{i,j}$
of the signed Entringer triangle satisfy relations
$$
\tilde\sharp_{i,j}=\tilde\sharp_{i,j-1}+\tilde\sharp_{i-1,j-1} \tag1.2
$$
with boundary conditions $\tilde\sharp_{i,1}=0$ for $i=2k$, $\tilde\sharp_{i,i}=0$
for $i=2k+1$,
$k>0$, $\tilde\sharp_{1,1}=1$. General triangles satisfying relation (1.2) with
arbitrary boundary conditions
were first studied more than 120 years ago by Seidel \cite{Se}. In particular, he
proved that
the ratio of exponential generating functions for the numbers on the right and on
the left sides of such
a triangle equals $e^x$. More recently such triangles where studied, from the
combinatorial point of
view, in \cite{DV, Du1, Du2}. In particular, it is proved in \cite{DV} that the
exponential generating function
for a Seidel triangle is equal to $e^yF(x+y)$, where $F(x)$ is the corresponding
function for the left side of the triangle.

The case of general up-down sequences was addressed by de Bruijn \cite{dB}
(see also \cite{Vi} for another version
of the same result). Let $\sharp_{i,j}^\P$ be the number of permutations of length $i$
whose last entry equals $j$ and whose up-down sequence equals $\P=(p_1,p_2,\dots)$.
He proved that these
numbers satisfy the following equations:
$$\alignedat 3
&\sharp_{i,j}^\P=\sharp_{i,j-1}^\P+\sharp_{i-1,j-1}^\P, &\quad &\sharp_{i,1}^\P=0,&
\qquad &\text{if \;$p_{i-1}=1$}\\
&\sharp_{i,j}^\P=\sharp_{i,j+1}^\P+\sharp_{i-1,j}^\P, &\quad &\sharp_{i,i}^\P=0,&
\qquad &\text{if \;$p_{i-1}=0$}.
\endalignedat
$$
with $\sharp_{1,1}^\P=1$. Evidently, for $\P=(10)^*$ one gets the Entringer
relations (1.1). As before, these
equations can be represented graphically as a triangle, and the direction in which
one has to advance along the rows
of the triangle is governed by the sequence $\P$. We call this triangle the
{\it de Bruijn triangle\/} corresponding
to the up-down sequence $\P$. A de Bruijn triangle is said to be {\it periodic\/}
if the corresponding up-down
sequence is periodic.

Let $\P$ be a periodic up-down sequence with period $p$ of length
$m>1$, and let $i_1<i_2<\dots< i_r$ be the locations of zeros in
$p$. Without loss of generality we assume that $i_r=m$ (otherwise
we consider instead of $\P$ the up-down sequence $\bar\P=(\bar
p_i)^*$, where $\bar p_i=1-p_i$ for $1\ls i\ls m$; evidently, the
de Bruijn triangle for $\bar\P$ is obtained from that for $\P$ by
the reflection in the vertical axis).

The {\it signed de Bruijn triangle\/} is obtained from the ordinary
de Bruijn triangle by multiplying  its $i$th row by
$$
\e_i=(-1)^{\bar p_1+\bar p_2+\dots+\bar p_{i-1}}, \qquad i\gs 1.
\tag1.3
$$
The corresponding exponential generating function is defined by
$$
F^\P(x,y)=\sum_{i\gs1}\sum_{j=1}^i \e_i\sharp_{ij}^\P
\frac{x^{i-j}y^{j-1}}{(i-j)!(j-1)!}. \tag1.4
$$

\proclaim{Theorem 1} The exponential generating function of the signed periodic de
Bruijn triangle corresponding
to the up-down sequence $\P$  is given by  $F^\P(x,y)=e^yf^\P(x+y)$, where
$$
f^\P(t)=\frac{\det \bar M^\P(t)}{\det M^\P(t)}
$$
and $M^\P(t)$ and $\bar M^\P(t)$ are $r\times r$ matrices
$$
M^\P(t)=\pmatrix \fy_{m,0}&\fy_{m,m-i_{1}}&\fy_{m,m-i_{2}}&\ldots&\fy_{m,m-i_{r-1}}\\
\fy_{m,i_{1}}&\fy_{m,0}&\fy_{m,m+i_{1}-i_2}&\ldots&\fy_{m,m+i_1-i_{r-1}}\\
\fy_{m,i_{2}}&\fy_{m,i_2-i_{1}}&\fy_{m,0}&\ldots&\fy_{m,m+i_2-i_{r-1}}\\
\vdots&\vdots&\vdots&\ddots&\vdots\\
\fy_{m,i_{r-1}}&\fy_{m,i_{r-1}-i_1}&
\fy_{m,i_{r-1}-i_2}&\ldots&\fy_{m,0}
\endpmatrix $$
and
$$\bar M^\P(t)=\pmatrix
1            & 1 & 1 &\ldots& 1\\
\fy_{m,i_{1}}&\fy_{m,0}&\fy_{m,m+i_{1}-i_2}&\ldots&\fy_{m,m+i_1-i_{r-1}}\\
\fy_{m,i_{2}}&\fy_{m,i_2-i_{1}}&\fy_{m,0}&\ldots&\fy_{m,m+i_2-i_{r-1}}\\
\vdots&\vdots&\vdots&\ddots&\vdots\\
\fy_{m,i_{r-1}}&\fy_{m,i_{r-1}-i_1}&
\fy_{m,i_{r-1}-i_2}&\ldots&\fy_{m,0}
\endpmatrix $$
with $\fy_{m,j}=\fy_{m,j}(t)$.
\endproclaim

In particular, for signed Entringer numbers one has $m=2$, $r=1$, and hence
$f^{\{10\}^*}(t)=\fy_{2,0}^{-1}(t)=1/\cosh t$, thus recovering the Arnold formula
for $A(x,y)$. Moreover,
the same techniques allows to obtain generating functions for other Seidel
triangles with periodic
boundary conditions, such as the triangle for Genocchi numbers (see \cite{DV}).
It can be also extended to
pairs of Seidel triangles with periodic boundary conditions, such as Arnold
triangles $L(\beta)$ and $R(\beta)$
for Springer numbers (see \cite{Ar2, Du2}), thus recovering several combinatorial
results obtained in \cite{Sp, Ar2}; see \S 2 for details.

As an immediate corollary of Theorem~1 we get the above mentioned
results of \cite{CGJN} and \cite{St2}  concerning the generating
functions for the number of permutations with a given up-down
sequence.

\proclaim{Corollary 1} Let $r$ be even, then the exponential
generating function for the numbers $\sharp_n^\P$ is equal
$$
1+\frac{\det \widetilde {M}^\P(t)}{\det M^\P(t)},
$$
where $\widetilde {M}^\P(t)$ is an $r\times r$ matrix
$$
\widetilde {M}^\P(t)=\pmatrix
\psi_{m,0}&\psi_{m,1}&\psi_{m,2}&\ldots&\psi_{m,r-1}\\
\fy_{m,i_{1}}&\fy_{m,0}&\fy_{m,m+i_{1}-i_2}&\ldots&\fy_{m,m+i_1-i_{r-1}}\\
\fy_{m,i_{2}}&\fy_{m,i_2-i_{1}}&\fy_{m,0}&\ldots&\fy_{m,m+i_2-i_{r-1}}\\
\vdots&\vdots&\vdots&\ddots&\vdots\\
\fy_{m,i_{r-1}}&\fy_{m,i_{r-1}-i_1}&
\fy_{m,i_{r-1}-i_2}&\ldots&\fy_{m,0}
\endpmatrix $$
with
$$
\psi_{m,j}=\e_{i_j+1}+\sum_{i\ne i_j}\e_{i+1}\fy_{m,m+i-i_j}
$$
and $\fy_{m,j}=\fy_{m,j}(t)$.
\endproclaim

To get a similar result for odd $r$, it suffices to consider $\P$
as a periodic sequence with a period of length $2m$.

 Let us now consider the asymptotic behavior of the numbers
$\sharp^\P_{i,j}$. It was observed without a proof in \cite{Ar2,
p.~18} that the even rows of the ordinary Entringer triangle
approximate, after an appropriate normalization, the function
$\sin x$ on the interval $[0,\pi/2]$, while the odd rows
approximate $\cos x$. Exact statements with the first two
correction terms can be found in \cite{SY}.

We generalize this result to arbitrary periodic de Bruijn
triangles.

\proclaim{Theorem 2} For any $l$, $0\ls l\ls m-1$, one has
$$
\lim_{n\to\infty, \frac{j}{mn+l}\to t}
\frac{\sharp_{mn+l,j}^\P\lambda^{mn+l}}{(mn+l-1)!}=c_{m,l}u_l^\P(t),
$$
where $c_{m,l}$ is a constant depending only on $m$ and $l$, and
$u_l^\P$ is the normalized first eigenfunction  of the two-point
spectral problem
$$
u^{(m)}=(-1)^r\lambda^mu
$$
with $m$ homogeneous  boundary conditions
$$\alignat2
&u^{(i)}(0)=0&\qquad &\text{if \;$p_{m+l-i-1}=1$},\\
&u^{(i)}(1)=0& \qquad &\text{if \;$p_{m+l-i-1}=0$},
\endalignat
$$
where $r$ is the number of zeros in the period $p$.

The first eigenvalue of the above spectral problem is the minimal
absolute value among the solutions of the equation
$$
\det  M^\P(t)=0,
$$
and the eigenfunction $u_l^\P$ is normalized by the condition
$\max_{0\ls t\ls 1}u_l^\P(t)=1$.
\endproclaim

In particular,  for the Entringer numbers $\sharp_{i,j}$ one gets
the sine law of \cite {SY}:
$$\align
&\lim_{k\to \infty,\; \frac {j}{2k+1}\to t}
\frac{\sharp_{2k+1,j}\left(\frac\pi2\right)^{2k+1}}{(2k)!} =
c_{2,1}\cos{\frac{\pi t}{2}},\\
&\lim_{k\to \infty,\; \frac {j}{2k}\to t}
\frac{\sharp_{2k,j}\left(\frac\pi2\right)^{2k}}{(2k-1)!} =
c_{2,0}\sin{\frac{\pi t}{2}};\endalign
$$
it is shown in \cite{SY} that $c_{2,0}=c_{2,1}=2$.

The starting point of this research was the result by the first
and the third author that the numbers $\sharp_{i,j}^\P$ for
$\P=(1^20^2)^*$  arise naturally in counting real rational
functions of a certain type, see \cite{SV}. The authors are
grateful to Max--Planck--Institut f\"ur Mathematik, Bonn for
     its hospitality in September 2000 and to the Royal Institute of
     Technology, Stockholm for the financial support of the visit of
     A.~V. to Stockholm in July-August 2001. Sincere thanks goes to
R.~Ehrenborg, A.~Laptev, H.~Shapiro, P.~Yuditski and A.Volberg for
useful discussions.

\heading {\S 2. Proofs} \endheading

\demo{Proof of Theorem 1}
Define $\tilde\sharp_{i,j}^\P=\e_i\sharp_{i,j}^\P$. Since by (1.3), $\e_i$ equals
$\e_{i-1}$ if
$p_{i-1}=1$ and $-\e_{i-1}$ otherwise, we immediately get that $\tilde\sharp_{i,j}^\P=
\tilde\sharp_{i,j-1}^\P+\tilde\sharp_{i-1,j-1}^\P$. In terms of the generating function
$F^\P$ defined by (1.4) this relation translates to $F^\P=\partial F^\P/\partial y-
\partial F^\P/\partial x$. The general solution of this equation is given by
$F^\P(x,y)=e^yf^\P(x+y)$, where $f^\P$ is a function of one variable to be defined
from the boundary conditions. Let $F^\P_L$ and $F^\P_R$ be the restrictions of
$F^\P$ to the left and the right side of the signed de Bruijn triangle. It follows
from the
above discussion that $F^\P_L(t)=f^\P(t)$ and $F^\P_R(t)=e^tf^\P(t)$. Evidently,
there exist unique representations
$$
f^\P(t)=\sum_{i=0}^{m-1}t^if^\P_i(t^m), \qquad
e^tf^\P(t)=\sum_{i=0}^{m-1}t^ig^\P_i(t^m);
$$
besides, $e^t=\sum_{i=0}^{m-1}\fy_{m,i}(t)$. Combining these expressions together
we get
$$\multline
\pmatrix g^\P_0(t^m)\\
        tg^\P_1(t^m)\\
         t^2g^\P_2(t^m)\\
         \vdots      \\
        t^{m-1}g^\P_{m-1}(t^m)\endpmatrix=\\
\pmatrix \fy_{m,0} & \fy_{m,m-1} & \fy_{m,m-2} & \dots & \fy_{m,1}\\
\fy_{m,1} & \fy_{m,0} & \fy_{m,m-1} & \dots & \fy_{m,2}\\
\fy_{m,2} & \fy_{m,1} & \fy_{m,0} &  \dots & \fy_{m,3}\\
\vdots & \vdots & \vdots & \ddots & \vdots\\
\fy_{m,m-1} & \fy_{m,m-2} & \fy_{m,m-3} & \dots &
\fy_{m,0}\endpmatrix
\pmatrix f^\P_0(t^m)\\
        tf^\P_1(t^m)\\
         t^2f^\P_2(t^m)\\
         \vdots      \\
        t^{m-1}f^\P_{m-1}(t^m)\endpmatrix,\endmultline\tag2.1
$$
where $\fy_{m,i}=\fy_{m,i}(t)$.
Besides, relations $p_{i}=0$ for $i=i_1,\dots, i_{r-1}$ imply
$\tilde\sharp^\P_{km+i+1, km+i+1}=0$ for $k=0,1,\dots$, and hence
$g^\P_{i}=0$ for $i=i_1,\dots,i_{r-1}$. Similarly, $p_m=0$ together
with $\tilde\sharp^\P_{1,1}=1$ imply $g^\P_{0}=1$,
and  $p_{i}=1$ for $i\ne i_1,\dots, i_{r}$ imply
$\tilde\sharp^\P_{km+i+1, 1}=0$ for $k=0,1,\dots$, and hence
$f^\P_{i}=0$ for $i\ne0,i_1,\dots,i_{r-1}$.
Therefore (2.1) is reduced to
$$
\pmatrix 1\\
        0\\
         0\\
         \vdots      \\
        0\endpmatrix=
\pmatrix \fy_{m,0}&\fy_{m,m-i_{1}}&\fy_{m,m-i_{2}}&\ldots&\fy_{m,m-i_{r-1}}\\
\fy_{m,i_{1}}&\fy_{m,0}&\fy_{m,m+i_{1}-i_2}&\ldots&\fy_{m,m+i_1-i_{r-1}}\\
\fy_{m,i_{2}}&\fy_{m,i_2-i_{1}}&\fy_{m,0}&\ldots&\fy_{m,m+i_2-i_{r-1}}\\
\vdots&\vdots&\vdots&\ddots&\vdots\\
\fy_{m,i_{r-1}}&\fy_{m,i_{r-1}-i_1}&
\fy_{m,i_{r-1}-i_2}&\ldots&\fy_{m,0}
\endpmatrix
\pmatrix f^\P_0(t^m)\\
        t^{i_1}f^\P_{i_1}(t^m)\\
         t^{i_2}f^\P_{i_2}(t^m)\\
         \vdots      \\
        t^{i_{r-1}}f^\P_{i_{r-1}}(t^m)\endpmatrix,
$$
and the result follows.
\qed
\enddemo

The Seidel triangle $\{\gamma_{i,j}\}$ for signed Genocchi numbers
presented in \cite{DV} is satisfies periodic boundary conditions
$\gamma_{2k,2k}=0$ for $k=1,2,\dots$, $\gamma_{2,2}=1$ and
$\gamma_{2k+1,1}+\gamma_{2k+1,2k+1}=0$ for $k=0,1,\dots$. Thus we
arrive to equations
$$
\pmatrix -f_0(t^2)\\
           t\endpmatrix=
\pmatrix \fy_{2,0}(t) & \fy_{2,1}(t)\\
         \fy_{2,1}(t) & \fy_{2,0}(t)\endpmatrix
\pmatrix f_0(t^2)\\
         tf_1(t^2)\endpmatrix,
$$
which mean that the generating function for this triangle is $2(x+y)e^y/(e^{x+y}+1)$.

In a slightly different way one can treat signed versions of
Arnold's pairs of triangles $\{L(\beta), R(\beta)\}$, $\{L(b),
R(b)\}$, and  $\{L(d), R(d)\}$ involving Euler and Springer
numbers (see \cite{Ar2, Du2}). The first pair consists of
triangles $\{\beta^L_{i,j}\}$ and $\{\beta^R_{i,j}\}$ satisfying
periodic boundary conditions
$$\align
&\beta^L_{2k+1,1}=0,\quad k=1,2,\dots, \qquad \beta^L_{1,1}=1,\\
&\beta^L_{k,k}=\beta^R_{k,1}, \quad k=0,1,\dots, \\
&\beta^R_{2k,2k}=0,\quad k=1,2,\dots.
\endalign
$$
Their signed versions obtained by multiplying the $i$th row by
$(-1)^{(i-1)(i-2)/2}$ are Seidel triangles.

Similarly, the second
 pair consists of triangles
$\{b^L_{i,j}\}$ and $\{b^R_{i,j}\}$ satisfying periodic boundary
conditions
$$\align
&b^L_{2k+1,1}=0,\quad k=0,1,\dots, \\
&b^L_{k,k}=b^R_{k,1}, \quad k=2,3,\dots,\qquad b^R_{1,1}=1, \\
&b^R_{2k,2k}=0,\quad k=1,2,\dots.
\endalign
$$
To get a pair of Seidel triangles one has to multiply the $i$th
row by $(-1)^{i(i-1)/2}$.

Finally, the third
 pair consists of triangles
$\{d^L_{i,j}\}$ and $\{d^R_{i,j}\}$ satisfying periodic boundary
conditions
$$\align
&d^L_{2k,1}=0,\quad k=1,2,\dots, \\
&d^L_{k,k}=d^R_{k,1}, \quad k=2,3,\dots,\\
&d^R_{2k+1,2k+1}=0,\quad k=1,2,\dots, \qquad d^R_{1,1}=1.
\endalign
$$
To get a pair of Seidel triangles one has to multiply the $i$th
row by $(-1)^{(i-1)(i-2)/2}$.

Let $f^\beta(t)=f^\beta_0(t^2)+tf^\beta_1(t^2)$ and
$g^\beta(t)=g^\beta_0(t^2)+tg^\beta_1(t^2)$ be the restrictions of the generating
functions for the signed versions of
$L(\beta)$ and $R(\beta)$ to the left and right sides, respectively.
Then we get the following equations:
$$
\pmatrix g^\beta_0(t^2)\\
           0\endpmatrix=
\pmatrix \fy_{2,0}(t) & \fy_{2,1}(t)\\
         \fy_{2,1}(t) & \fy_{2,0}(t)\endpmatrix^2
\pmatrix 1\\
         tf^\beta_1(t^2)\endpmatrix,
$$
and hence the generating functions for the triangles $L(\beta)$ and $R(\beta)$ are
$$
F^\beta_L(x,y)=\frac{e^{-2x-y}}{\cosh 2(x+y)},\qquad
F^\beta_R(x,y)=\frac{e^{-x}}{\cosh 2(x+y)}.
$$

For similar restrictions $f^b(t)$ and $g^b(t)$ one has
$$
\pmatrix g^b_0(t^2)\\
           0\endpmatrix=
\pmatrix \fy_{2,0}(t) & \fy_{2,1}(t)\\
         \fy_{2,1}(t) & \fy_{2,0}(t)\endpmatrix^2
\pmatrix 1\\
         tf^b_1(t^2)\endpmatrix+
\pmatrix \cosh t\\
         t^{-1}\sinh t\endpmatrix,
$$
and hence the generating functions for the for the signed versions of
triangles $L(b)$ and $R(b)$ are
$$
F^b_L(x,y)=-\frac{e^y\sinh (x+y)}{\cosh 2(x+y)},\qquad
F^b_R(x,y)=\frac{e^{-x}\cosh (x+y)}{\cosh 2(x+y)}.
$$

Finally, for $f^d(t)$ and $g^d(t)$ one has
$$
\pmatrix            1\\
           tg^d_1(t^2)\endpmatrix=
\pmatrix \fy_{2,0}(t) & \fy_{2,1}(t)\\
         \fy_{2,1}(t) & \fy_{2,0}(t)\endpmatrix^2
\pmatrix f^d_0(t^2)\\
                  0\endpmatrix+
\pmatrix \cosh t\\
         \sinh t\endpmatrix,
$$
and hence the generating functions for the for the signed versions of
triangles $L(d)$ and $R(d)$ are
$$
F^d_L(x,y)=\frac{e^y(1-\cosh (x+y))}{\cosh 2(x+y)},\qquad
F^d_R(x,y)=\frac{e^{x+2y}(1-\cosh (x+y))}{\cosh 2(x+y)}.
$$

\demo{Proof of Corollary 1} It is enough to observe that for $r$
even, the exponential generating function for $\sharp_n^\P$ is
equal to $1+\sum_{i=0}^{m-1}\e_{i+1}t^i(f^\P_i+g^\P_i)$, and to
evaluate this sum according to the proof of Theorem~1.
\qed
\enddemo

\demo{Proof of Theorem~2} Consider the space $PC[0,1]$ of
piecewise continuous functions on the segment $[0,1]$ with the
norm $||f||=\max_{0\ls t\ls 1}|f(t)|$. Define two operators, $S^0$
and $S^1$, taking $PC[0,1]$ to itself:
$$
S^0f(x)=\int_0^x f(t) dt,\qquad S^1f(x)=\int_x^1 f(t) dt.
$$
It is easy to see that both $S^0$ and $S^1$ can be written as
integral operators
$$
S^if(x)=\int_0^1K^i(x,y) f(y) dy,\quad  i=0,1,
$$
with the kernels
$$
K^0(x,y)=\cases 1, \quad x\gs y,\\
                0, \quad x<y, \endcases
$$
and $K^1(x,y)=1-K^0(x,y)$. Besides, consider two families of
integral operators $S^i_n$, $i=0,1$, with the kernels $K^i_n(x,y)$
given by
$$
K^0_n(x,y)=\cases 1,\quad \lfloor(n+1)x\rfloor\gs\lfloor
ny\rfloor,\\
0,\quad\text{otherwise,}\endcases
$$
and $K^1_n(x,y)=1-K^0_n(x,y)$. It is easy to see that
$$
||S^i-S^i_n||\ls\frac cn \tag 2.2
$$
for some constant $c>0$.

Consider the operator $S^\P=S^{\bar p_m} S^{\bar p_{m-1}}\cdots
S^{\bar p_1}$. Evidently, $S^\P$ is compact, as a product of
compact operators. To study the spectral properties of $S^\P$, we
make use of the infinite-dimensional Perron-Frobenius theory, as
presented in \cite{KR}.

\proclaim{Lemma~1} For any $f\in PC[0,1]$,
$$
S^\P f=\mu\psi(f)u^\P+Af, \tag 2.3
$$
where $\mu>0$, $\psi$ is a strictly positive functional, $S^\P
u^\P=\mu u^\P$, $||u^\P||=1$, $Au^\P=0$, $\psi(Af)=0$ and
$\lim_{n\to\infty}||A^n||^{1/n}<\mu$.
\endproclaim

\demo{Proof} Recall that we have assumed without loss of
generality that $p_m=0$. Let $s$ be the largest index satisfying
relations $p_{m+1-s'}=0$ for $1\ls s'\ls s$.
 Define a linear operator $R\:PC[0,1]\to PC[0,1]$ by
$Rf(x)=(1-x)^sf(x)$, and let $\R$ be the image of $PC[0,1]$ under
$R$. We consider $\R$ as a Banach space in the induced norm, that
is, $||f||_\R=\max_{0\ls x\ls 1}\frac{f(x)}{(1-x)^s}$.

Let $K$ be the cone of nonnegative functions in $\R$, that is,
$f\in\R$ belongs to $K$ if $f(x)\gs0$ for $x\in[0,1]$. Evidently,
the closure of the linear span of $K$ coincides with $\R$. We say
that a linear operator $A\: \R\to\R$ is {\it strongly positive\/}
if it preserves $K$ and takes any boundary point of $K$ to an
inner point. It is easy to see that $f\in K$ is an inner point of
$K$ if $\lim\inf_{x\to1}f(x)/(1-x)^s>0$.

Consider the restriction of $S^\P$ to $\R$, which we denote by
$S^\P_\R$. Evidently, $S^\P_\R$ is compact and strongly positive.
Hence, by Theorem~6.3 of \cite{KR}, the maximal eigenvalue $\mu$
of $S^\P_\R$ is positive, the corresponding eigenspace is
one-dimensional, and all the other eigenvalues lie strictly inside
the circle of radius $\mu$. Since $S^\P$ takes the whole $PC[0,1]$
to $\R$, the same $\mu$ is the maximal eigenvalue of $S^\P$, and
all the above properties of $S^\P_\R$ remain valid for $S^\P$. It
follows from the proof of Theorem~6.3 in \cite{KR} that (2.3)
holds. \qed\enddemo

By differentiating equation $S^\P u^\P=\mu u^\P$ $m$ times we see
that $u^\P$ satisfies equation $\mu u^{(m)}=(-1)^ru$ with the
boundary conditions $u^{(i)}(0)=0$ if $p_{m-i}=0$ and
$u^{(i)}(1)=0$ if $p_{m-i}=1$. Put $\lambda=\mu^{-1/m}$; then
$\lambda$ is the minimal positive solution of the equation $\det
\widehat{M}^\P(\xi t)=0$ with $\xi^m=(-1)^r$ and
$$
\widehat{M}^\P(t)=\pmatrix \fy_{m,0}&\fy_{m,m-i_{r-1}}&\fy_{m,m-i_{r-2}}&\ldots&\fy_{m,m-i_{1}}\\
\fy_{m,i_{r-1}}&\fy_{m,0}&\fy_{m,m+i_{r-1}-i_{r-2}}&\ldots&\fy_{m,i_{r-1}-i_1}\\
\fy_{m,i_{r-2}}&\fy_{m,m+i_{r-2}-i_{r-1}}&\fy_{m,0}&\ldots&\fy_{m,i_{r-2}-i_1}\\
\vdots&\vdots&\vdots&\ddots&\vdots\\
\fy_{m,i_{1}}&\fy_{m,m+i_1-i_{r-1}}&
\fy_{m,m+i_1-i_{r-2}}&\ldots&\fy_{m,0}
\endpmatrix $$
with $\fy_{m,j}=\fy_{m,j}(t)$. It is easy to see that $M^\P$ can
be obtained from $\widehat{M}^\P$ by the transformation
$(k,l)\mapsto (r+2-k,r+2-l)$, and hence $\lambda$ is the minimum
absolute value of the solutions of the equation $\det M^P(t)=0$.

Observe that the condition on the norms of $||A^n||$ in (2.3)
means that there exists $k$ such that $||A^k||<\mu^k$. Denote
$\alpha=\frac1{\mu^k}||A^k||<1$ and $T=(\frac1\mu S^\P)^k$.

We are approximating the operator $S^\P$ with the help of
operators $S^\P_n$ defined by
$$
S^\P_{n+1}=S^{\bar p_m}_{(n+1)m}S^{\bar p_{m-1}}_{(n+1)m-1}\cdots
S^{\bar p_1}_{ nm+1}.
$$
It follows immediately from (2.2) that
$$
||S^\P-S^\P_n||\ls \frac {c'}{n}+O\left(\frac1{n^2}\right).\tag2.4
$$
Finally, to approximate $T$ define operators
$T_{n+1}=\frac1{\mu^k}S^\P_{(n+1)k}S^\P_{(n+1)k-1}\cdots
S^\P_{kn+1}$; it follows from the above inequality that
$$
||T_n-T||\ls\frac {c''}n +O\left(\frac1{n^2}\right).
$$

For any function $f$ define a sequence $\{f_n\}$ by
$f_n=T_nf_{n-1}$ with $f_0=f$.

\proclaim{Lemma~2} Let $\{||f_n||\}$ be bounded, then  the
sequence $f_n=f_n/||f_n||$ converges to  $u^\P$ as $n\to\infty$.
\endproclaim

\demo{Proof} By (2.3), each of $f_i$ can be uniquely represented
as $f_i=f_i^u+f_i^A$, where $f_i^u$ is a multiple of $u^\P$ and
$f_i^A$ belongs to the maximal subspace invariant under $A$ and
not containing $u^\P$. Therefore,
$$
f_n=T_nf_{n-1}=(T+(T_n-T))(f_{n-1}^u+f_{n-1}^A)= f_{n-1}^u+
A^kf_{n-1}^A+ (T_n-T)(f_{n-1}^u+f_{n-1}^A),
$$
which together with $||u^\P||=1$ gives
$$\gather
||f_n^u||\gs ||f_{n-1}^u||- ||(T_n-T)f_{n-1}||,\\
||f_n^A||\ls||A^kf_{n-1}^A||+||f_{n-1}^u||\cdot||(T_n-T)u^\P||+||(T_n-T)f_{n-1}^A||.
\endgather
$$
Recall that $\{||f_n||\}$ is bounded, and hence
$$
||f_n^u||\gs
||f_{n-1}^u||-\left(\frac{c'''}n+O\left(\frac1{n^2}\right)\right)\gs
\beta||f_{n-1}^u||,
$$
where $\beta<1$ can be chosen arbitrary close to $1$ for $n$ big
enough. Therefore,
$$\multline
\delta_n=\frac{||f_n^A||}{||f_n^u||}\ls \frac
{||f_{n-1}^u||}{\beta||f_{n-1}^u||}\left(\frac{c''}{n}+
O\left(\frac{1}{n^2}\right)\right)+\frac{||f_{n-1}^A||}{\beta||f_{n-1}^u||}
\left(\alpha+\frac{c''}{n}+O\left(\frac1{n^2}\right)\right)\ls\\
\frac{c''}{\beta n} +\left(\frac{\alpha}{\beta}+\frac{c''}{\beta
n}\right)\delta_{n-1}+O\left(\frac1{n^2}\right).
\endmultline
$$
Therefore, either $\delta_{n-1}\ls n^{-1/2}$, or
$\delta_{n-1}>n^{-1/2}$ and
$$
\delta_n\ls\left(\frac{c''}{\beta
\sqrt{n}}+\frac{\alpha}{\beta}+\frac{c''}{\beta
n}\right)\delta_{n-1}+O\left(\frac1{n^2}\right) \ls
\alpha'\delta_{n-1}
$$
for some constant $\alpha'<1$. In any case $\delta_n\to0$ as
$n\to\infty$, and  hence the sequence $\{f_n/||f_n^u||\}$
converges to a multiple of $u^\P$, which implies the convergence
of $\{f_n/||f_n||\}$ to $u^\P$. \qed\enddemo

 For an arbitrary finite
sequence $a=\{a_1,\dots,a_k\}$ define a piecewise constant
function $\hat a\in PC[0,1]$ whose value equals $a_i$ on the
interval $\left[\frac{i-1}k,\frac ik\right)$ for $i=1,\dots,k-1$
and $a_k$ on $\left[\frac{k-1}k,1\right]$. Let
$\sharp_k=\{\sharp_{k,1},\dots,\sharp_{k,k}\}$. It is easy to see
that $\hat\sharp_{k+1}=kS^{\bar p_k}_k\hat\sharp_k$, which means
that
$$
S_{n+1}^\P\hat\sharp_{mn+1}=\frac{(mn)!}{(m(n+1))!}\hat\sharp_{m(n+1)+1}.\tag2.5
$$
Observe that the sequence of functions
$g_n=\frac{\hat\sharp_{mn+1}\lambda^{mn+1}}{(mn)!}$ is bounded.
Indeed, $||\hat\sharp_{mn+1}||=\sharp_{mn+1,1}$. The exponential
generating function for the numbers $\sharp_{mn+1,1}$ is
calculated in Corollary~1. Since the numerator of the
corresponding expression is a polynomial in Olivier functions,
which converge in the whole plain, the numbers
$\frac{\sharp_{mn+1,1}}{(mn)!}$ grow asymptotically as
$\frac{\gamma}{\lambda^{mn+1}}$, hence for $n$ big enough one has
$\frac{\sharp_{mn+1,1}\lambda^{mn+1}}{(mn)!}<\gamma'$, and
therefore the sequence $\{g_n\}$ is bounded. Moreover, (2.5) can
be rewritten as $\lambda^mS^\P_{n+1}g_{mn+1}=g_{m(n+1)+1}$.
Therefore, Lemma~2 applies, and $g_{nm+1}\to u^\P$ as
$n\to\infty$.

Combining the above results we get
$$
\lim_{n\to\infty, \frac{j}{mn+1}\to t}
\frac{\sharp_{mn+1,j}^\P\lambda^{mn+1}}{(mn)!}=c_{m,1}u_1^\P(t),
$$
with $c_{m,1}=\gamma'$, and hence Theorem~2 is proved for $l=1$.

To get the proof for the other values of $l$ one has to consider,
instead of $S^\P$, a different operator: $S^{\bar
p_{m+l-1}}S^{\bar p_{m+l-2}}\dots S^{\bar p_l}$. Its properties
are identical to those of $S^\P$; to prove this one has to use
operators $R$ and $L\: f(x)\mapsto x^sf(x)$, depending on the
value of $p_{m+l-1}$. \qed
\enddemo

\enddocument